\input amstex
\documentstyle{amsppt}

\loadeufb
\loadeusb
\loadeufm
\loadeurb
\loadeusm

\magnification =\magstep 1
\refstyle{1}
\NoRunningHeads

\topmatter
\title On uniformization of compact Kahler spaces 
\endtitle

\author Robert  Treger \endauthor
\address Princeton, NJ 08540  \endaddress
\email roberttreger117{\@}gmail.com \endemail
\endtopmatter

\document 
\head
1.  Introduction
\endhead 

The aim of the present note is to extend author's uniformization theorem in \cite{29}
 to compact  Kahler spaces with mild singularities and  establish a kind of
 rigidity of  
 their universal coverings. We study compact   Kahler spaces with large, 
residually finite and nonamenable fundamental groups and mild singularities 
by passing to their universal coverings  where there are plenty of bounded
harmonic functions (see Lyons and Sullivan\cite{19} and Toledo \cite{25}). 

\medskip

(1.1) Throughout the note $X$ will be a connected, compact, normal, Cohen - Macaulay Kahler space of dimension $n$. If $X$ is, in addition,  projective we then assume the canonical class $\eusm K_X$ is a $\bold Q$-divisor. Also, we assume the fundamental group $\Gamma:=\pi_1(X)$ is large, residually finite and nonamenable. 
\medskip
(1.2) Let $U_X$ denote the universal covering of $X$.  Recall that  $\pi_1(X)$ is {\it large}\/ if and only if $U_X$ contains no proper holomorphic subsets of positive dimension \cite{17}.

Recall that a countable group $G$ is called {\it amenable}\/ if there is on $G$ a finite additive, translation invariant nonnegative probability measure (defined for all subsets of $G$). Otherwise, $G$ is called {\it nonamenable}.  

In 1970s, it became clear  that in the classification 
of projective manifolds of dimension at least three, one has to consider
varieties with mild singularities.
Let $Z$ denote a normal holomorphic space.
We say that $Z$ has at most 1-rational singularities if
$Z$ has a resolution of singularities $\psi: Y\rightarrow Z$ such that
$R^1 \psi_*\Cal O_Y =0$. It is well known that $Z$ has {\it rational singularities}\/ if and only if  $R^i \psi_*\Cal O_Y =0 \; (\forall i\geq 1)$; equivalently, $Z$ is Cohen - Macaulay and $\psi_*\omega_Y = \omega_Z$, where $\omega$ denotes the dualizing sheaf. 
Set  $Z^{\roman {reg}}: = Z\backslash Z_{\roman{sing}}$.

Bounded domains in $\bold C^n$ have the (classical) Bergman metric.
As was already observed by Bochner \cite{2} (also see Kobayashi
 \cite{15}, \cite{16}), one can view the Bergman metric in a bounded domain $D$ via a natural  embedding of $D$ into a Fubini space $\bold F_\bold C(\infty,1)$ (infinite-dimensional projective  space).
Furthermore, the image of $D$ does not intersect a hyperplane at infinity. 
Taking, then,  the Calabi flattening out of $\bold F_\bold C(\infty,1)$ into  $\bold F_\bold C(\infty,0)$ (infinite-dimensional Hilbert space) \cite{3, Chap.\;4}, our $D$ will be holomorphically embedded into $\bold F_\bold C(\infty,0)$.

Now, we say that $U_X$ has a {\it generalized}\/ Bergman metric (Bergman metric in a singular space) if we get a natural  embedding of $U_X$ into a Fubini space $\bold F_\bold C(\infty,1)$ 
given by functions-sections of a pluricanonical bundle.  Furthermore, its image does not
 belong to a proper subspace of $\bold F_\bold C(\infty,1)$ and does not intersect a hyperplane at infinity.

The infinite-dimensional complex projective spaces $\bold F_\bold C(\infty,1)$ as well as their real counterparts $\bold F_\bold R(2\infty,1)$ play a prominent 
role in the present note. Often, $U_X$ can be flatten out into $\bold F_\bold C(r,0)$ with $r<\infty$.

Let $X$ and $Z$ satisfy the assumptions (1.1) and are homeomorphic:
$X\approx  Z$. We assume they have generalized Bergman metrics and can be flatten
out into $\bold F_\bold C(r_X,0)$ and   $\bold F_\bold C(r_Z,0)$, respectively,
with $r_X, r_Z<\infty$ (e.g., $X$ and $Z$ are nonsingular \cite{27}).
Then we say  $U_X$ is {\it rigid}\/ if the induced mapping $U_X \rightarrow U_Z$ is homotopic to a biholomorphism.

The aim of the present note is to prove the following theorem.

\proclaim{Theorem (Uniformization)}  We assume, in addition to the assumptions (1.1), that $X$ has at most rational singularities. Then  
\roster
\item $U_X$ is Stein,
\item  $\eusm K^m_X$ is  very ample for an integer $m\gg0$, and $U_X$ has a generalized Bergman metric,

\item $U_X$ is rigid.
\endroster
\endproclaim

First, we establish  that $U_X$ is Stein provided $X$ is a projective variety. 
Second, we will show that our Kahler space $X$ is a Moishezon space.
Therefore $X$ is a projective variety by a theorem of Namikawa \cite{21} because $X$ has at most 1-rational singularities (a generalization of  Moishezon's theorem).
We, then, establish (2) and (3).
\medskip

\head
2.  Preliminaries
\endhead

(2.1) {\it Kahler spaces}.\/ Recall that Kahler spaces were defined by  Grauert \cite{10, Sect.\;3.3}: Let  $Y$ be an arbitrary reduced holomorphic space with a Hermitian metric $ds^2 =\sum g_{i k} dz_i d\bar z_k$ in $Y\backslash Y_{\roman{sing}}$ (with continuous coefficients). Then it is called a Kahler metric in $Y$ if and only if 
for every point $a\in Y$, there is an open neighborhood $\Cal Y$ of $a$ and a strictly plurisubharmonic function  $p$ on 
$\Cal Y$ such that $g_{ik} = {\partial^2 p \over\partial z_i \partial \bar z_{ k}}$ in $\Cal Y\backslash \Cal Y_{\roman{sing}}$. 

A function $u:
Y \rightarrow \bold R \cup \{-\infty \}$ is said to be weakly plurisubharmonic (equivalently, plurisubharmonic by Fornaess and Narasimhan \cite{7}) if for any holomorphic map $\varepsilon : \Delta\rightarrow Y$ of the unit disk into $Y$, its pullback $\varepsilon ^*u$ is either subharmonic or equal to $-\infty$.

Later, Moishezon \cite{20} defined a Kahler space as follows: An arbitrary reduced
holomorphic space $Y$ is said to have a  Kahler metric if the exists and open covering  $\{\Cal Y_\mu\}$ of $Y$ and a system of strictly plurisubharmonic functions 
$p_\mu $ of class $C^2$, with each $p_\mu$ defined on $\Cal Y_\mu$, such that 
$p_\mu-p_\nu$ is pluriharmonic function on $\Cal Y_\mu \cap \Cal Y_\nu$. 
Both definition coincide when $Y$ is normal. 

Thus, a normal holomorphic space $Z$ has a Kahler metric (in the sense of Grauert or Moishezon) if there is a Kahler metric $h$ on $Z^{\roman {reg}}$ , and every singular point has a neighborhood $\Cal Z \subset Z$ and a closed embedding $\Cal Z \subset \Cal U$, 
where $\Cal U$ is an open subset of 
an affine space, such that there is a Kahler metric $h^\prime$ on $\Cal U$
with $h^\prime|_{\Cal Z^{\roman {reg}}} = h$.
Note that Grauert and Moishezon considered not necessary real-analytic Kahler metrics.

If we assume, in addition, that $X$ is a projective variety then $U_X$ will be equipped with $\pi_1(X)$-equivariant Kahler metric  as in \cite{27, Sect.\;3}; see Lemma 3.1 below.

Recall the following well-known fact (see \cite{9, Lemma 1}). Let $(\Cal Y, a)$ be a normal isolated singularity with $\dim(\Cal Y, a) \geq 2$.  Let $f$ be a strictly plurisubharmonic function on $\Cal Y\backslash  a$. Then, for any sufficiently small neighborhood $\Cal V$ of $a$, there exists a strictly plurisubharmonic function $f^\prime$ on $\Cal Y$ such that $f=f^\prime$ on $\Cal Y\backslash \Cal V$.
\medskip

(2.2) {\it Di\'astasis}.\/ For complex {\it manifolds}\/ $M$ with a real-analytic Kahler metric, the diastasis was introduced by Calabi \cite{3, Chap.\;2} (see also
\cite{30, Appendix} as well as a brief review in \cite{27}). 
This notion can be generalized to  normal Cohen - Macaulay holomorphic spaces $M$ with real-analytic Kahler metric. 
If $\dim M=1$ then we will consider (not necessary normal) reduced holomorphic 
spaces with real-analytic Kahler metric (in the sense of Moishezon).

 We will mention only the {\it fundamental property of the diastasis},\/ namely, it is  inductive on normal Cohen - Macaulay holomorphic subspaces with real-analytic Kahler metric (see a proof for manifolds in \cite{3, Chap.\;2, Prop.\;6} as well as \cite{30, Lemma 1.1}). 
Clearly, the other properties of the diastasis mentioned in \cite{27, (2.2), (2.2.1),
(2.3.2), (2.3.2.1), (2.3.2.2)} can be generalized to real-analytic Kahler spaces. 
\smallskip

(2.2.1) {\it Example}. Let  
 $H$ denote a complex Hilbert space, i.e. $\bold F_\bold C(N,0)$, where
$1\leq N \leq \infty$. Then $\bold F_\bold C(N,0) \cong \bold F_\bold R(2N,0)$.
The diastasis of $\bold F_\bold C(N,0)$ with respect to the reference point  $\bold Q$ is 
$D(\bold Q,p)= \sum^N_{i=1} |z_i(p)|^2$, It is equal to the diastasis of
$\bold F_\bold R(2N,0)$, i.e., $\sum^N_{i=1} (|x_i(p)|^2+|y_i(p)|^2)$, where 
$z_i=x_i+\sqrt{-1}y_i$. 
\smallskip

Now, let $\bold P_\bold C(H)$ be a complex projective space,
i.e. $\bold F_\bold C(N,1)$. In the canonical coordinates around $\bold Q$, the diastasis of  $\bold F_\bold C(N,1)$ is
$$
D(\bold Q,p) = \log (1+ \sum^N_{i=1} |z_i(p)|^2).
$$
The diastasis of the real projective space $\bold F_\bold R(2N,1)$ equals 
the diastasis of $\bold F_\bold C(N,1)$. We can consider Calabi's flattening out
of  $\bold F_\bold C(N,1)$ as well as  $\bold F_\bold R(2N,1)$.
\smallskip

In the present note,  we often employ induction on dimension. Recall the following well-known local Bertini theorem: Let $a\in M$ be a  point of an arbitrary normal
Cohen  -  Macaulay holomorphic space $M$.
If $\dim M\geq 3$ then the general local (or global in the quasi-projective case) hyperplane section  through $a$ is also normal and Cohen - Macaulay.
\medskip

(2.3) {\it Harmonic mappings.}\/ For the early applications of harmonic mappings to the study of mapping between Kahler manifolds, see the surveys by Siu \cite{23} and Yau
\cite{31} and references therein.

Also, Anderson \cite{1} and Sullivan \cite{24} constructed non-constant bounded harmonic functions on any simply-connected complete Riemannian manifold whose sectional curvature is bounded from above by a negative number.

Harmonic mappings between admissible Riemannian polyhedra were considered by Eells - Fuglede \cite{6} and in  papers by Serbinowski, Fuglede
and others (see \cite{8}  and references therein) where the concept of
 energy is defined in  
the spirit of Korevaar - Schoen \cite{18}.  Also see the paper by
 Gromov - Schoen \cite{14}, 
and the survey articles by Toledo \cite{26} and Delzant - Gromov \cite{5}.

 Recall that an arbitrary normal holomorphic space is an admissible Riemannian polyhedron. 
It has been established the existence and uniqueness of a solution to the variational Dirichlet problem for harmonic mappings into the target 
space with nonpositive curvature \cite{8}. The latter can be applied to rigidity questions.

\medskip

(2.4) {\it Reduction theory}\/ (also see \cite{10, Sect.\;2}).  In 19th Century, Riemann and Poincar\'e studied compact complex manifolds of dimension one with non-Abelian fundamental group by passing to the universal covering of the manifold. In dimension one, simply-connected manifolds are very simple. However, only in 1907, Koebe and
Poincar\'e established, independently, that the universal coverings are discs $\Delta\subset \bold C$ provided the genus of the manifold is at least two.

If the dimension of the manifold is at least two then its universal covering can be very complicated.  Furthermore, many compact Kahler manifolds can be simply connected.

So, one approaches the classification of higher-dimensional  algebraic manifolds
along the lines developed by the Italian school of algebraic geometry in the late
19th Century and the early 20th Century. 
In 1930s, Hodge employed  transcendental methods in his study of algebraic varieties as well as more general Kahler manifolds   bypassing the uniformization problem as well. His approach was proceeded by the Lefschetz  topological methods in algebraic geometry.

In the 19th Century, the  blowing downs were studied 
by Cremona, Del Pezzo and others, and were later applied to the classification of algebraic surfaces by Castelnuovo and Enriques. It is still one of the main tools  in the classification of compact complex algebraic manifolds.

In 1956, Remmert 
proved that a holomorphically convex holomorphic space can be transformed into a  holomorphically complete holomorphic space. 
Precisely, given a holomorphically convex holomorphic space $Y$, under some  additional assumptions, Remmert proved  that there exists a holomorphically complete $Z$ and proper holomorphic mapping 
$$
\Phi : Y \longrightarrow Z
$$ 
such that  if $\Cal U\subset Z$ is any open set, $\Cal V=\Phi^{-1}(\Cal U)$, and $f$ is a holomorphic function on $\Cal V$ then
there exist a holomorphic function  $g$ on $\Cal U$ such that  $f=g\circ \Phi$. 

The paire $(\Phi,Z)$ is called the holomorphic reduction of $Y$. If $Y$ is normal then $Z$ is also normal.
In 1960, Cartan proved that Remmert's additional assumptions 
can be removed. Moreover, for any $a\in Z$, $\Phi^{-1}(a)$ is a connected holomorphic space. 
\medskip

(2.5) {\it Exceptional sets}\/ \cite{10, Sect.\;2.4, Def.\;3}. Nowhere discrete compact holomorphic subset $\eufm X$ of a holomorphic space $Y$ is called {\it exceptional}\/  if there is a holomorphic space $Z$ and a proper holomorphic mapping $\Phi : Y\rightarrow Z$ transforming $\eufm X$ onto a discrete set $D\subset Y$, mapping biholomorphically the set $Y\backslash \eufm X$  onto $Z\backslash D$, and, in addition, for every neighborhood $\Cal U =
 \Cal U(D) \subset Z$ and every holomorphic function $g$ on $\Cal V: =\Phi^{-1}(\Cal U)$
there exists a holomorphic function $f$ on $\Cal U$ such that $g=f\circ \Phi$.

Grauert proved that $\eufm X \subset Y$ is exceptional if and only if 
there exists a strongly pseudoconvex 
neighborhood $\Cal U=\Cal U( \eufm X) \Subset Y$ such that its maximal compact 
holomorphic subset equals $\eufm X$ \cite{10, Sect.\;2.4, Theorem 5}.

According to Grauert, a vector bundle $\bold V$ over $X$ is said to be {\it weakly 
negative}\/ if its zero section $\eufm X$ has a strongly pseudoconvex neighborhood $\Cal U
=\Cal U (\eufm X)
\Subset \bold V$ \cite{10, Sect.\;3.1, Def.\;1}. He also established that $\bold V$ is
weakly negative if and only if its zero section is exceptional (see \cite{10, Sect.\;3.1, Theorem 1}). Also, Grauert proved the following extension to holomorphic spaces of the celebrated Kodaira embedding theorem: Given  a compact holomorphic space $Y$  with
 a weakly negative line bundle $\Cal L$ over $Y$, then  $Y$ is a projective variety and
$\Cal L^{-1}$ is an ample line bundle over $Y$ \cite{10, Sect.\;3.2, Theorem 2
and its proof}.
\medskip

(2.6) {\it Shafarevich's conjecture.}\/ In 1960s, H.\;Wu conjectured that the universal 
covering of a compact Kahler manifold with negative sectional curvature is a bounded 
domain in an affine space. He had also established that the covering is Stein.
The Wu conjecture was recently established by the author 
in \cite{28} as well as in \cite{29}, provided
the fundamental group is residually finite.

In 1971, Griffiths considered the following problem \cite{11, Question 8.8}: Let $D$
be a complex manifold and $\Gamma \subset {\roman {Aut}}(D)$ a properly discontinuous group of automorphisms such that $D/\Gamma$ is a quasi-projective algebraic variety. Then do the meromorphic functions separate points of $D$? Is $D$
meromorphically convex?

In the early 1970s, based on Remmert's reduction, Shafarevich conjectured that the universal covering of  any projective manifold is holomorphically convex \cite{22}.  According to Shafarevich, if the conjecture were established then, roughly speaking, one could reduce  the study of the universal covering of a projective manifold to a compact case and a holomorphically complete case. 

In the early 1990s, Campana (Kahler case) and Koll\'ar (projective case) constructed independently  a meromorphic map
$$
Y --> Sh(Y) =\Gamma(Y)
$$
such that for a general point $a\in Y$, the fiber $Y_a$ passing through $a$ is the largest among the connected holomorphic subsets $\Cal A$ of $Y$ containing
the point $a$ such that 
the natural map $i_* : \pi_1 (\Cal A^{nor}) \rightarrow \pi_1(Y)
$ has a finite image, where $i_*$ is induced by the inclusion of $\Cal A$ in $Y$
composed with the normalization map $\Cal A^{nor} \rightarrow \Cal A$.

This map is called the Shafarevich map by Koll\'ar and the $\Gamma$-reduction of $Y$ by Campana. They also established the relationship of this map and the corresponding map for the universal covering $U_Y$ with (unsolved at the time) Shafarevich's conjecture.
The Shafarevich conjecture, established by 
the author a few years ago provided the fundamental group is residually finite, 
implies the meromorphic  map is regular and $a\in Y$ can be any point.

\medskip

\head
3. Projective case
\endhead

  We will establish the following version of the Shafarevich conjecture for singular projective varieties. We need the following lemma.

\proclaim{Lemma 3.1 (Extension of metric)} Let $Z$ be an arbitrary  normal Cohen - 
Macaulay holomorphic space  of dimension $n\geq 1$. Given an arbitrary real-analytic Kahler metric on $ Z^{\roman{reg}}$, we get a unique real-analytic Kahler metric on $Z$ whose restriction on $ Z^{\roman{reg}}$ is the given one. 
\endproclaim

\demo{Proof of Lemma 3.1} We will prove the lemma by induction on dimension of $Z$. 
If $Z=Z^{\roman{reg}}$ there is nothing to prove.
First, we consider the case $\dim Z=2$. Then $Z$ has only isolated singularities. 
As in Section 2.1, we get a unique Kahler metric on $Z$. 

This Kahler metric will be real analytic. Indeed, in a small 
 neighborhood $\Cal V$ of a point $q\in Z$, the potentials of this metric  produce a unique holomorphic function in a neighborhood of the diagonal of
$\Cal V^{\roman{reg}} \times \bar\Cal V^{\roman{reg}}$ as follows. This holomorphic
function arises from the  diastasic potential of $\Cal V^{\roman{reg}}$  which 
is a real-analytic function over the whole $\Cal V^{\roman{reg}}$.
Next, the holomorphic function can be
extended to a holomorphic function  on the diagonal of $\Cal V \times \bar\Cal V$. The latter holomorphic function produces a real-analytic function on $\Cal V$ by  realification (i.e., the inverse of  comlexification) which  is a potential of our Kahler metric at the point $q$.

If $\dim Z\geq 3$ we apply the local Bertini theorem and induction on dimension.  
Let $p\in Z$ be a singular point. By induction hypothesis, we get a real-analytic
potential (diastasic potential) in a small neighborhood inside 
 a general local hyperplane section $v_p$ through $p\in Z$. We obtain a holomorphic 
function in a neighborhood of $p\times \bar p$ (in $v_p \times v_{\bar p}$) employing complexifications.

By Hartogs' theorem (separate analyticity\! $\implies$\! joint analyticity), we obtain holomorphic functions in the  points of  diagonal of $Z \times \bar Z$.
Finally, we produce real-analytic potentials in every point of $Z$ by realifications 
because the  obtained real-analytic functions will be strictly plurisubharmonic 
by Fornaess - Narasimhan \cite{7}.

This proves the lemma.
\enddemo

\proclaim{Proposition 1} Let $\phi : Y\hookrightarrow \bold P^r$ be a connected,  normal, Cohen - Macaulay  projective variety of dimension $n\geq 1$ with 
large and residually finite $\pi_1(Y)$. If the genus of a general curvilinear section of $Y$ is at least two then $U_Y$ is a Stein space.
\endproclaim

\demo{Proof of Proposition 1} We assume the map $\phi$ is given by a very ample
bundle $\Cal L$. As in \cite{27, Sect.\;3}, we obtain the real-analytic $\pi_1(Y)$-invariant
 Kahler metric
$\Lambda_\Cal  L$ on $U_{Y^{\roman {reg}}}$. By Lemma 3.1, we extend this metric to $U_Y$.

As in \cite{27, Appendix, Lemma A}, we establish the prolongation of the diastasic potential of our metric at a point of  $U_Y$ over the whole $U_Y$.

Finally, as in \cite{27, Appendix, (A.2)}, we conclude the proof of holomorphic completeness by reduction to the one-dimensional case because one-dimensional non-compact holomorphic spaces are Stein.

\enddemo
\medskip

\head
4. Kahler case
\endhead

  We will establish the following Kahler version of Proposition 1.

\proclaim{Proposition 1$^\prime$} Let $Y$  be a connected, compact, normal, Cohen - Macaulay,   Kahler space of dimension $n\geq 1$ with large and residually finite $\pi_1(Y)$. 
If $\pi_1(Y)$ is nonamenable and $Y$ has at most 1-rational singularities then $U_Y$ is a 
Stein space and $Y$ is projective.
\endproclaim

\demo{Proof of Proposition 1$^\prime$} As in \cite{29}, we will establish that that $Y$ is projective by showing that $Y$ is a Moishezon space.  
Let $Har(U_Y)$ ($Har^b(U_Y)$  be the space of  harmonic functions  (bounded  harmonic functions) on $U_Y$. The space $Har^b(U_Y)$ contains
non-constant harmonic functions by Lyons 
and Sullivan \cite{19}. 
In fact, the space $Har^b(U_Y)$ is infinite dimensional by Toledo \cite{25}.
Their proof is stated for Riemannian manifolds but it works in our case. One can consider the transition density of Brownian motion  on $U_{Y^{\roman {reg}}}$ or $U_Y$.

We will integrate  the {\it bounded} harmonic functions  with respect to the measure
$$
dv:=p_{U_Y}(s,x,\bold Q) d\mu = p_{U_Y}(x) d\mu,
$$
where $\bold Q \in U_Y$ is a fixed point (a so-called center of $U_Y$), $d\mu$ is the   Riemannian measure and $p_{U_Y}(x):= p_{U_Y}(s,x,\bold Q)$ is the 
heat kernel. We obtain a pre-Hilbert space of bounded harmonic 
 functions (compare \cite{27, Sect.\;2.4 and Sect.\;4}); note that all 
 bounded harmonic  functions are square integrable, i.e. in $L_2(dv)$.

 We observe that the latter pre-Hilbert  space has a completion 
in the Hilbert space $H$ of all harmonic  $L_2(dv)$ functions:
$$
H := \biggl\{h \in Har(U_Y)\;\biggl |\; \parallel h\parallel^2_H:=  \int_{U_Y}
|h (x)|^2 dv =\int_{U_Y}|h (x)|^2 p_{U_Y}(x) d\mu < \infty   \biggl\}.
$$ 
Let $H^b \subseteq H$ be the Hilbert subspace generated by $Har^b(U_Y)$.
These Hilbert spaces are separable infinite dimensional and have  reproducing kernels.
 
Let $\{\phi_j\}\subset Har^b(U_Y)$ be an orthonormal basis of $H^b$.
We obtain a continuous, even smooth,  finite $\Gamma$-energy $\Gamma$-equivariant mapping into $(H^b)^*$:
$$
g: U_Y \longrightarrow (H^b)^*\quad (\subset \bold P_\bold R (H^b)^*)),
\qquad u \mapsto (\phi_0(u), \phi_1(u), \dots). \tag{1}
$$
The group $\Gamma$ acts  isometrically  on 
$(H^b)^*$ via 
$
\psi \mapsto (\psi\circ\gamma).
$

We assume $g$ is harmonic; otherwise, we replace $g$ by a harmonic mapping
homotopic to $g$. Clearly, the proposition is valid when $\dim Y =1$ even if $Y$ is
not necessary normal. So we assume $\dim Y \geq 2$ in the following lemmas.

\enddemo

\proclaim{Lemma 4.1} With assumptions of Proposition 1$^{\prime}$, $g$ will produce a pluriharmonic mapping $g^{fl}$. There exists a natural holomorphic mapping $g^{h} : U_Y\longrightarrow \bold F_\bold C(\infty,0)$.
\endproclaim

\demo{Proof of Lemma 4.1}  We define a harmonic  $\Gamma$-equivariant mapping as
follows
$$
g^{fl}: =\eurb S_{g(\bold Q)}\!\circ \! g: U_Y \longrightarrow (H^b)^*.
$$
 We have applied the mapping $g$ followed by 
 the Calabi {\it flattening out}\/ $\eurb S_{g(\bold Q)}$ 
of the real projective space $\bold F_\bold R
(2\infty,1)$ from  ${g(\bold Q)}$ 
into the Hilbert space \cite{3, Chap.\;4, p.\;17}.

By  \cite{3, Chap.\;4, Cor.\;1, p.\;20}, the whole  $\bold F_\bold R(2\infty,1)$,
except the antipolar hyperplane $A$ of $g(\bold Q)$, can be flatten out into $\bold F_\bold R (2\infty,0)$. The image of $g$ does not intersect the antipolar hyperplane $A$ of $g(\bold Q)$. Thus we have introduced a flat metric in a large (i.e.\;outside $A$)  neighborhood  of  ${g(\bold Q)}$ in $ \bold P_\bold R((H^b)^*)$.

Since the mapping $g^{fl}$ has finite $\Gamma$-energy, it is pluriharmonic; this  is a  special case of the well-known theorem of Siu. 
Since $U_Y$ is simply connected, we obtain the natural holomorphic mapping $g^{h}:
U_Y \longrightarrow \bold F_\bold C (\infty,0)$.

\enddemo

\proclaim{Lemma 4.2} Construction of a complex line bundle $\Cal L_Y$ on $Y$ and 
its pullback on $U_Y$, denoted by $\Cal L$.
\endproclaim

\demo{Proof of Lemma 4.2} We take a  point $u\in U_Y$. Let $v :=
g^h(u) \in \bold F_\bold C(\infty, 1)$, where $ \bold F_\bold C(\infty, 1)$ is the complex projective space (see (1)). We consider the linear system of hyperplanes in
$\bold F_\bold C(\infty, 1)$ through $v$ and its proper transform on $U_Y$. 
We consider only the moving part. The projection
on $Y$ of the latter linear system on $U_Y$ will produce a linear system on $Y$. 

A connected component of
a {\it general} member of the latter linear system on $Y$ will be an irreducible divisor 
$D$ on $Y$ by Bertini's theorem. The corresponding line bundle will be 
$\Cal L_Y:=\Cal O_Y(D)$ on $Y$. 
\enddemo

\proclaim{Lemma 4.3} Conclusion of the proof of Proposition 1$^{\prime}$ by induction on\/ $\dim Y$.
\endproclaim

\demo{Proof of Lemma 4.3}  
By the Campana-Deligne theorem \cite{17, Theorem 2.14}, $\pi_1(D)$ will be nonamenable. 
We proceed by induction on $\dim Y$. As above, let $\dim Y \geq 2$.
Let $q=q(n)$ be an appropriate integer.

We get a  global holomorphic function-section  
$f$ of $\Cal L^q$ corresponding to a bounded pluriharmonic function (see Lemma 4.1 
and \cite{27, Sect.\;4}).

We will define a $\Gamma$-invariant Hermitian quasi-metric  on sections of $\Cal L^q$  below.
Furthermore,  $f$  is $\ell^2$ on orbits of $\Gamma$, and it  is not  identically zero on any orbit because, otherwise, we could have
replaced
$U_Y$ by $U_Y\backslash B$, where the closed holomorphic subset $B\subset U_Y$ is the union of those orbits on which $f$ had vanished \cite{17, Theorem 13.2, Proof of Theorem 13.9}.

One  can show that  $f$ satisfies the above conditions by taking linear systems of curvilinear sections of $U_Y$ through  $u\in U_Y$ and their projections on $Y$ 
(see the
proof of Lemma 4. 2 above), since  the statements are trivial in dimension one.
The required Hermitian quasi-metric on $\Cal L^q_{\roman {reg}}:= \Cal L^q|U_{Y^{\roman {reg}}}$ is also defined by 
induction on dimension with the help of the Poincar\'e residue map \cite{13,
pp.\;147-148}.

The condition $\ell^2$ on orbits of $\Gamma$ is a local property on $Y^{\roman {reg}}$. We
get only a Hermitian quasi-metric on $\Cal L^q_{\roman {reg}}$ (instead of a Hermitian metric). Precisely, we get Hermitian metrics over small neighborhoods of
points of $Y^{\roman {reg}}$,   and on the intersections of neighborhoods, they will 
 differ by constant multiples (see \cite{17, Chap.\;5.13}).

For $\forall k>N\gg 0$, the Poincar\'e series are holomorphic sections over $U_{Y^{\roman {reg}}}$:
$$
P(f^k)(u):= \sum_{\gamma\in \Gamma} \gamma^*f^k(\gamma u) \qquad (u\in U_{Y^{\roman {reg}}})
$$
and they do not vanish for infinitely many $k$ (see \cite{17, Sect.\;13.1, Theorem 13.2}).

Finally, we can apply Gromov's theorem, precisely, its generalization by Koll\'ar 
(see \cite{12, Corollary 3.2.B, Remark 3.2.B$'$} and 
\cite{17, Theorem 13.8, Corollary 13.8.2, Theorem 13.9, Theorem 13.10}).
So, $Y$ is a Moishezon space hence a projective variety.

The Lemma 4.3 and Proposition 1$^\prime$ are established.

\enddemo
\medskip

\head
5. Conclusion of the proof of theorem
\endhead

\proclaim{Proposition 2} With assumptions of the theorem, $\eusm K^m_X$ is ample for an integer $m >0$, and $U_X$ has a generalized Bergman metric.
\endproclaim

\demo{Proof of Proposition 2}  For simplicity, we assume that $\eusm K_X$ is a line bundle.

It  follows from the Kodaira
embedding theorem that the sheaf $\eusm K_X$ is ample, provided $X$ is nonsingular
(see  \cite{27, Sect.\;4}). 

In the general case, the ampleness of $\eusm K_X$  will follow from Grauert's extension of the Kodaira 
theorem (see Section 2.5 and \cite{27, Sections 4, 5} with $\Cal L:=\eusm K_X^{-1}$).

Thus, we have to show $\eusm K_X^{-1}$ is weakly negative, i.e., its zero section 
$\eufm X $ is exceptional. Namely, we have to prove the existence of a strongly pseudoconvex neighborhood $\Cal U=\Cal U(\eufm X) 
\Subset \eusm K_X^{-1}$ 
such that its maximal compact holomorphic subset equals $\eufm X$. 

We apply Grauert's generalization of Kodaira's result 
(see \cite{10, Sect.\;3, Proof of Theorem 1 (negativity $\implies$ exceptional}). 
We will be able to choose a covering $ \{\Cal V_\tau \}$ of $X$ and strictly plurisubharmonic 
functions $-\log p_\tau$ on each $\Cal V_\tau$ so that the following Hermitian form
is strictly negative on $X$:
$$
\sum  {\partial^2\log p_\tau (z,\bar z) \over \partial z_\alpha \partial \bar z_\beta}
dz_\alpha dz_\beta.  
$$
 The latter form is negative of the Ricci form of the volume form 
on $X$ obtained as follows. 
First, we consider the resolution of singularities of $X$: $X^\prime \rightarrow X$ (and the corresponding pullback $U^\prime_X \rightarrow U_X$). 
As in Section 4, employing bounded harmonic functions, we consider the 
$\Gamma$-invariant volume form on $U^\prime_X$ (compare \cite{27, Sect.\;4}). 

We get a $\Gamma$-invariant  Hermitian metric on
$\eusm K^{-1}_{U_X^\prime}$. We, then, get   Hermitian metrics on 
$\eusm K_{X^\prime}^{-1}$ and $\eusm K_{X}^{-1}$ and obtain the desired
$p_\tau$'s.
Hence $\eufm X$ is an exceptional section.

Now, the existence of a generalized Bergman metric on $U_X$ is established as in 
\cite{27, Sections 5.1-5.3} because, for a {fixed}\/ large integer $q$, $\eusm K^q$ is very ample on every finite Galois covering of $X$. If $X$ is nonsingular the very ampleness follows from \cite{4}, 
and, in our case, we resolve the singularities of $X$ and apply the latter assertion.

\enddemo

\proclaim{Proposition 3} $U_X$ is rigid.
\endproclaim

\demo{Proof of Proposition 3} Let $X$ and $Z$ satisfy the assumptions of the theorem and
$X\approx Z$. We get a $\Gamma$-invariant topological equivalence $U_X \approx U_Z$.
Both $U_X$ and $U_Z$ are equipped with the  $\Gamma$-invariant  generalized Bergman metrics and can be flatten out into the compact spaces $\bold F^X_\bold C (r_X,0)$
and $\bold F^Z_\bold C (r_Z,0)$, respectively.
We get the embeddings into the corresponding Fubini spaces
$$
g_X: U_X \hookrightarrow \bold F^X_\bold R (2r_X,0) \quad {\text{and}} \quad
g_Z: U_Z \hookrightarrow \bold F^Z_\bold R (2r_Z,0).
$$
 
We consider the towers of finite Galois coverings: 
$$
X\subset \cdots \subset X_i \subset \cdots \subset  U_X \qquad 
{\text{and}}\qquad 
Z\subset \cdots \subset Z_i \subset \cdots \subset  U_Z;
$$
$$ 
\bigcap_i Gal(U_X/X_i) =  \bigcap_i Gal(U_Z/Z_i) = \{1\}. 
$$
We have assumed the topological equivalence $X_i \approx Z_i$. For the 
corresponding fundamental domains
$F_{X_i} \subset U_X$ and $F_{Z_i} \subset U_Z$, we have the topological equivalences $F_{X_i} \approx F_{Z_i}$ for all $i$. 

Now, we can and will assume the {\it embeddings}
$$
F_{Z_i} \hookrightarrow U_Z \qquad (\forall i), 
$$
where $U_Z$ and the interior of each $F_{Z_i}$ are considered with the 
above generalized Bergman metric, and $F_{Z_i}$ is the Dirichlet fundamental domain
\cite{27, (2.3.3.2)}.

From existence and uniqueness of a solution to the variational Dirichlet problem
 for harmonic mappings, we get the  mappings
$$
F_{X_i} \longrightarrow F_{Z_i} \quad (\forall i)
$$
that are harmonic in the interiors of $F_{X_i}$ and coincide with the given continuous 
mappings on the boundaries. These mappings will be pluriharmonic in the interiors.

 Thus we get a  pluriharmonic mapping
$$
U_X \longrightarrow U_Z \subset
\bold F^Z_\bold R (2r_Z,0).
$$

Similarly, we get a pluriharmonic mapping
$$
U_Z \longrightarrow U_X \subset
\bold F^X_\bold R (2r_X,0).
$$

 Since $U_Y$ and $U_Z$ are simply connected, the above pluriharmonic mappings
produce holomorphic mappings, and $U_X$ is biholomorphic to $U_Z$ by 
uniqueness of a solution  to  the variational Dirichlet problem:
$$
\bold F^X_\bold C (r_X,0) \supset U_X \rightleftharpoons U_Z \subset
\bold F^Z_\bold C (r_Z,0).
$$

This proves Proposition 3 and the Theorem.
\enddemo

\enddocument

Stefan.Kebekus@math.uni-freiburg.de
schoen@math.stanford.edu

hwbllmnn(at)mpim-bonn.mpg.de

•	Prof. Dr. Stefan Kebekus
Arbeitsgruppe: Algebraische Geometrie
Tel. +49 761 203-5536   Raum: 425 (Eckerstr. 1)
Stefan.Kebekus@math.uni-freiburg.de


This mapping  will be pluriharmonic. Thus we get a holomorphic mapping
$$
h^{YZ} : U_Y \rightarrow U_Z \subset \bold F_\bold C(r,0)
$$

for an appropriate integer $r<\infty$.

Similarly we get $h^{ZY} : U_Z \rightarrow U_Y$. 
It follows from the existence and uniqueness of the solution to the variational Dirichlet problem for harmonic mappings that $U_Y$ is biholomorphic to $U_Z$.


\ref
\key T1 \by  R. Treger \pages
\paper Uniformization
\yr
\jour arXiv:math.AG/1001.1951v4
 \endref

\ref 
 \key FK  \by  J. Faraut, S. Kaneyuki, A. Kor\'anyi, Q.-k. Lu, G. Roos
\book Analysis and geometry on complex
homogeneous domains
\publ Birkh\"auser, Boston
\yr 2000
\endref

\ref 
 \key FG  \by  K. Fritzsche, H. Grauert 
\book From Holomorphic Functions to Complex manifolds
\publ Graduate Texts in Mathematics Ser., Springer 
\yr 2010
\endref

\ref
\key Sha \by B. V. Shabat
\book  Introduction to Complex Analysis 
\publ Nauka, Moscow
\yr 1969
\lang Russian
\endref

 Both $Y$ and $U_Y$ are equipped with the Kahler metrics with curvature bounded from above by $\kappa <0$.
From the  topological equivalence $Y\simeq Z$, we get corresponding
homotopic $\Gamma$-equivariant harmonic 
mappings $U_Y\rightarrow U_Z$ and $U_Z\rightarrow U_Y$, respectively. 
As in Lemma 4.1, we get pluriharmonic $\Gamma$-equivariant embeddings into the corresponding Fubini spaces
$$
g_Y : U_Y \hookrightarrow \bold F_\bold R^Y (0,\infty), \quad g_Z : U_Z \hookrightarrow \bold F_\bold R^Z (0,\infty).
$$
We observe that $\Gamma$-spaces $\bold F_\bold R^Y (0,\infty)$ and 
$\bold F_\bold R^Z (0,\infty)$ are not necessary assumed to be isomorphic.

Recall that $\Gamma$ is residually finite. From the existence and uiquenes of solution of the variational Dirichlet problem, we obtain a  one-to-one harmonic  mappings
$$
g_Y(U_Y) \rightleftharpoons g_Z(U_Z).
$$
Indeed, the above mappings arise as a limit of the mappings of the  Dirichlet fundamental domains (see \cite{TMet, Sect.\;2.3.3.2}) 
of finite coverings of $Y$ in $U_Y$ and the 
Dirichlet fundamental domains of the corresponding finite coverings of $Z$ in $U_Z$.

As in Lemma 4.1, the above mappings will be pluriharmonic, and  we get the desired
biholomorphism
$$
g^h_Y(U_Y)   \rightleftharpoons    g^h_Z(U_Z).
$$

This proves Proposition 3 and the Theorem.
\enddemo
  \enddocument

This proves Proposition 3 and the Theorem.
\enddemo

Anderson M., The Dirichlet Problem at Infinity for Manifolds of Negative
Curvature, J. Diff. Geom. 18 (1983), pp. 701-721.

Sullivan D., The Dirichlet Problem at Infinity for a Negatively Curved Manifold,
J. Diff. Geom. 18 (1983), pp. 723-732.

\enddocument

By  Grauert-Remmert and  Fornaess-Narasimhan, it can be uniquely extended to a  plurisubharmonic function on  $Y$ 
provided $\dim Y = 2$. 

Second, we apply the local Bertini's theorem for normal varieties ( see \cite{O-s}) and the fundamental property of diastasis. 
We can reduce the proof of lemma to the 2-dimensional case. 

Finally, because the extensions are unique we get the strictly plurisubharmonic functions in small neighborhoods of every point of $Y$. These strictly plurisubharmonic functions
will be real analytic
\enddemo

(3.2) We will establish the following generalization of the Shafarevich conjecture to singular varieties.

\demo{Proof} Let $\Cal L= \Cal L_X$ denote the very ample line bundle defining the map $\phi$. As in \cite{TMet},we can construct a Kahler metric on $U_X$. First, we get 
a proposition-definition similar to \cite{Proposition-Definition 2}. Second, we get a real analytic Kahler metric on $U_{X^{\roman {reg}}}$. We, then,  apply Lemma 1 and get the desired Kahler metric.

Finally, exactly as in \cite{T, Lemma A in Appendix}, we establish the prolongation of the diastasic potential of our metric on $X$. We conclude the proof of the proposition as in 
\cite{T, Appendix, (A.2)}.

\enddemo

\enddocument

We proceed as in \cite{T-metr} (we employ the same notation). We assume  $U_X$ and all
$X_i$'s are equipped with the Kahler metric $\Lambda _\Cal  L$. Then we establish the prolongation (see \cite{Lemma A, Apendix, T-metr}). We conclude the argument as in 
\cite{(A.2), Apendix, T-metr}.

\enddocument

In the \lq\lq Basic Algebraic Geometry\rq\rq, Shafarevich conjectured that that universal covering $\tilde M$ of a (projective)  manifold  $M$ is holomorphically convex. Thus one could reduce the study of $\tilde M$ to a compact case and a holomorphically complete case. In the early 1990s, Campana \cite{Cam} (Kahler case) and Koll\'ar \cite{Kol}(projective case) independently proved a kind of reduction for $M$ and $\tilde M$ and explained its relationship to the Shafarevich conjecture.


\vfill
\eject
For complex {\it manifolds}\/ $M$ with a real analytic Kahler metric. Let $\Phi$ denote a real analytic potential of the metric defined in a small neighborhood  $\Cal V \subset M$. Let $z=(z_1,\dots,z_n)$ be a coordinate system in $\Cal V$ and $\bar z=(\bar z_1,\dots,\bar z_n)$ a coordinate system in its conjugate neighborhood $\bar\Cal V\subset \bar M$. 
Let $(\bold p,\bold p)$ be a point on the diagonal of $M \times {\bar M}$ such that the neighborhood $\Cal V \times {\bar \Cal V} \subset M \times {\bar M}$ contains the point.

There exists a unique {\it holomorphic}\/ function $F$ on an open neighborhood of $(\bold p,\bold p)$ such that $F_{(\bold p,\bold p)}= \Phi_\bold p$.
Here $\Phi_\bold p$ is the germ at $\bold p\in M$ of our real analytic function, and $F_{(\bold p,\bold p)}$ is the germ of the corresponding holomorphic function (complexification of $\Phi_\bold p$ \cite {C, Chap.\;2}, \cite{U, Appendix}). 

One considers the sheaf $\Cal A^\bold R_M$ of germs of real analytic functions on $M$, and the sheaf $\Cal A^\bold C_{M \times \bar M}$ of germs of complex holomorphic functions on $M \times \bar M$. For each $\bold p\in M$, we get a natural inclusion $\Cal A^\bold R_{M,\bold p} \hookrightarrow \Cal A^\bold C_{M \times \bar M,(\bold p,\bold p)}$, called a complexification. The above equality is understood in this sense.

  Now, let $p$ and $q$ be two {\it arbitrary}\/ points of $\Cal V$ with coordinates $z(p)$ and $z(q)$. Let $ F(z(p), \overline{z( q)})$ denote the  complex holomorphic function on $\Cal V \times \bar \Cal V$ obtained from from $\Phi$. 
The {\it functional element of diastasis}\/ is defined as follows \cite{C, (5)}:
$$
D_M(p ,q):=F(z(p), \overline{z( p)})\,+\,F(z(q), \overline{z( q)})\,-\,F(z(p), \overline{z( q)})\,-\,F(z(q), \overline{z(p)}). \tag{2.2.1}
$$
We get the germ $D_M(p,q)\in \Cal A^\bold C_{M \times \bar M,(p,q)}$, and $D_M(p,q)$ is uniquely determined  by the Kahler metric, symmetric in $p$ and $q$ and real valued  \cite{C, Prop.\;1,\,2}.
The {\it real}\/ analytic function generated by the above functional element is called the diastasis \cite{C, p.\;3}. 
The diastasis approximates the square of the geodesic distance in the {\it small}\/ \cite{C, p.\;4}.
For $\bold C^r$ with its unitary coordinates,
$
D_{\bold C^r}(p,q) = \sum^r_{i=1} |z_i(q) -z_i(p)|^2.
$

The {\it fundamental property of the diastasis}\/ is that it is inductive on complex submanifolds \cite{C, Chap.\;2, Prop.\;6}.  

Now, let $\bold Q\in M$ be a {\it fixed}\/ point, and $z=(z_1, \dots, z_n)$ a local coordinate system in a small neighborhood $\Cal V_\bold Q \subset M$ with origin at $\bold Q$. The real analytic function $\tilde \Phi_\bold Q (z(p), \overline{z(p)}):= D_M(\bold Q,p)$ on $\Cal V_\bold Q$ is called the {\it  diastasic potential at}\/ $\bold Q$ of the Kahler metric. It is strictly plurisubharmonic function in $p$ \cite{C, Chap.\;2, Prop.\;4}. 

The {\it prolongation}\/ over $M$ of  the germ of   $\tilde \Phi_\bold Q (z(p), \overline{z(p)})$ at $\bold Q$ is a {\it function}\/ $\bold P_M:=\bold P_{M,\bold Q} \in  H^0(\Cal A^\bold R_M, M)$ such that,  for every $u\in M$,  $\bold P_{M}(u)$ coincides  with  $D_M(\bold Q, u) $ meaning $D_M(\bold Q, u) $, initially defined  in a neighborhood of $\bold Q$, can be extended over the {\it whole} $M$. Moreover, the germ of $\bold P_M$ at $u$ is the  diastasic potential of our metric at $u$. 

The  prolongation over $M$ of the germ of  diastasic potential  is not always possible, e.g., there are no strictly plurisubharmonic functions on $\bold P^1$.  
Now, let  $\bold P^N$ be a projective space with the Fubini-Study metric ($1\leq N \leq \infty$). For $\bold Q \in\bold P^N $, we consider Bochner canonical coordinates $z_1, \dots, z_N$ with origin at  $\bold Q$ on the complement of a hyperplane at infinity. By Calabi \cite{Chap.\;4, (27)}, 
$$
	D_{\bold P^r}(\bold Q,p)= \log\bigl(1+\sum^N_{\sigma=1}|z_\sigma(p)|^2 \bigl).
$$
In the homogeneous coordinates $\xi_0,\dots,\xi_N  $, where $z_\sigma := {\xi_\sigma/\xi_0}$, we get
$$
D_{\bold P^N}(\bold Q,p)= \log {\sum^N_{\sigma=0} |\xi_\sigma(p)|^2 \over |\xi_0(p)|^2}.
$$

\enddocument

\demo{Proof of Theorem} By a  theorem of Moishezon \cite{7}, it will suffice to establish that $X$ is a Moishezon manifold. 
In \cite{4, Sect.\;3},  Gromov uses his notion of Kahler
hyperbolicity to obtain holomorphic $L_2$ forms on $U_X$ and prove that $X$ is 
Moishezon. A priori, we do not know if there are holomorphic $L_2$ forms on $U_X$.

Set $\Gamma: =\pi_1(X)$. Let $\Cal L$ be an arbitrary complex line bundle on $U_X$. We will consider a section $f\in H^0(\Cal L^q,U_X)$ which is not assumed
to be  $L_p$, where $p < \infty$.
As in Koll\'ar \cite{5, Chap\;13.1}, we will employ $\ell^p$ sections $f$  {\it on orbits} of $\Gamma$ in place of $L_p$ sections.
Of course, we need a natural $\Gamma$-invariant Hermitian {\it quasi}-metric on $\Cal 
L^q$ (see the definition in the proof of Lemma 3).

Given an arbitrary $\Gamma$-invariant Hermitian metric on $U_X$, we get the
induced Riemannian metric on $U_X$ with the volume form $d\mu$. Since $\Gamma$ is nonamenable, we get non-constant  bounded harmonic functions 
on $U_X$ by Lyons and Sullivan \cite{6}. 
Employing their theorem, Toledo \cite{8} has established that the space of bounded harmonic functions as well as the space generated by bounded positive harmonic functions are infinite dimensional (see \cite{9, Sect.\;2.6}). 
Given $r$ linearly  independent functions $g_1,\dots,g_r$ on
$U_X$,  clearly  there exist $r$ points $ u_1,\dots,u_r
 \in U_X$ such that the vectors  
$
\langle g_1(u_i),\dots,g_r(u_i)\rangle\; (1\leq i\leq r)
$ 
are linearly independent. 

Let $Har(U_X)$ ($Har^b(U_X)$) be the space of  harmonic functions  (bounded  harmonic functions, respectively) on $U_X$.

In place of the standard $L_2(d\mu)$ integration with the standard Riemannian measure $d\mu$ on $U_X$, we will integrate  the {\it bounded} harmonic functions  with respect to the measure 
$$
dv:=p^2_{U_X}(s,x,\bold Q) d\mu,
$$
where $\bold Q \in U_X$ is a fixed point  and $p_{U_X}(s,x,\bold Q)$ is the heat kernel. Because all  bounded harmonic  functions are square integrable, i.e. in $L_2(dv)$, we obtain the pre-Hilbert space of bounded harmonic  functions (compare \cite{9, Sect.\;2.4 and Sect.\;4}).
 We observe that the latter pre-Hilbert  space has a completion in the Hilbert space $H$ of all harmonic  $L_2(dv)$ functions:
$$
H := \biggl\{h \in Har(U_X)\;\biggl |\; \parallel h\parallel^2_H:=  \int_{U_X}
|h (x)|^2 dv < \infty   \biggl\}. 
$$
Let $H^b \subseteq H$ be the Hilbert subspace generated by $Har^b(U_X)$.
These Hilbert spaces are infinite dimensional and have  reproducing kernels. 
The group $\Gamma$ acts  isometrically on $H^b: \psi \mapsto B[(\psi\circ\gamma)\cdot  Jac_\gamma]\; (\gamma\in \Gamma)$, where $Jac_\gamma$ is the complex 
Jacobian determinant of $\gamma$ and $B$ is the corresponding Bergman projection
(compare \cite{9, Sect.\;4.6} and references therein).

Let $\{\phi_j\}\subset Har^b(U_X)$ be an orthonormal basis of $H^b$.
We obtain a continuous, even smooth,  finite $\Gamma$-energy $\Gamma$-equivariant mapping 
$$
g: U_X \longrightarrow \bold P((H^b)^*)\qquad \big(u \mapsto [\phi_0(u) : \phi_1(u) : \dots]\big)
$$
by assigning to each point $u\in U_X$ the hyperplane of  $H^b$ consisting of functions vanishing at $u$.  Also we get a mapping $g: U_X \longrightarrow (H^b)^*,\;u\mapsto \psi(u)\, (\forall \psi \in H^b)$. 


We assume $g$ is harmonic; otherwise, we replace $g$ by a harmonic mapping
homotopic to g.
Let $\bold F_\bold C(\infty,0)$ denote the complex flat Fubini space, i.e.  a complex Hilbert space. 

\enddemo

\proclaim{Lemma 1} With assumptions of the theorem, $g$ will produce a pluriharmonic mapping $g^{fl}$. There exists a natural holomorphic mapping $g^h : U_X \longrightarrow \bold F_\bold C(\infty,0)$.
\endproclaim

\demo{Proof of Lemma 1}  We define a harmonic  $\Gamma$-equivariant mapping
$$
g^{fl}: =\eurb S_{g(\bold Q)}\cdot  g: U_X \longrightarrow (H^b)^*.
$$
 We have applied the mapping $g$ followed by 
 the Calabi {\it flattening out}\/ $\eurb S_{g(\bold Q)}$ (a generalized stereographic projection from   $g(\bold Q)$) of the real projective space $\bold F_\bold R(\infty,1)$  
into the Hilbert space \cite{2, Chap.\;4, p.\;17}.
By  \cite{2, Chap.\;4, Cor.\;1, p.\;20}, the whole  $\bold F_\bold R(\infty,1)$,
except the antipolar hyperplane $A$ of  $g(\bold Q)$, can be flatten out into $\bold F_\bold R (\infty,0)$. The image of $g$ does not intersect the antipolar hyperplane $A$ of $g(\bold Q)$. Thus we have introduced a flat metric in a large (i.e.\;outside $A$)  neighborhood  of  $g(\bold Q)$ in $ \bold P((H^b)^*)$.

Since the mapping $g^{fl}$ has finite $\Gamma$-energy, it is pluriharmonic; this  is a  special case of a theorem of Siu (see, e.g., \cite{1}). 
Since $U_X$ is simply connected, we obtain the natural holomorphic mapping $g^h:
U_X \longrightarrow \bold F_\bold C (\infty,0)$.

\enddemo

\proclaim{Lemma 2} Construction of a complex line bundle $\Cal L_X$ on $X$ and 
its pullback on $U_X$, denoted by $\Cal L$.
\endproclaim

\demo{Proof of Lemma 2} We take a  point $u\in U_X$. Let $v :=
g^h(u) \in \bold F_\bold C(\infty, 1)$, where $ \bold F_\bold C(\infty, 1)$ is the complex projective space. We consider the linear system of hyperplanes in
$\bold F_\bold C(\infty, 1)$ through $v$ and its proper transform on $U_X$. The projection
on $X$ of the latter linear system on $U_X$ will produce a linear system on $X$. 

A connected component of
a {\it general} member of the latter linear system on $X$ will be an irreducible divisor 
$D$ on $X$ by Bertini's theorem. The corresponding line bundle will be the desired
$\Cal L_X:=\Cal O_X(D)$ on $X$. 
\enddemo

\proclaim{Lemma 3} Conclusion of the proof of theorem by induction on $\dim X$.
\endproclaim

\demo{Proof of Lemma 3} 
By the Campana-Deligne theorem \cite{5, Theorem 2.14}, $\pi_1(D)$ will be nonamenable. 
We proceed by induction on $\dim X$, the case
$\dim X=1$ being trivial. Let $q=q(n)$ be an appropriate integer.

We get a  global holomorphic function-section  
$f$ of $\Cal L^q$ corresponding to a bounded pluriharmonic function (see Lemma 1 
and \cite{9, Sect.\;4}).
We will define a $\Gamma$-invariant Hermitian quasi-metric  on sections of $\Cal L^q$  below.
Furthermore,  $f$  is $\ell^2$ on orbits of $\Gamma$, and it  is not  identically zero on any orbit because, otherwise, we could have replaced
$U_X$ by $U_X\backslash B$, where the closed analytic subset $B\subset U_X$ is the union of those orbits on which $f$ had vanished \cite{5, Theorem 13.2, Proof of
Theorem 13.9}.

One  can show that  $f$ satisfies the above conditions by taking linear systems of curvilinear sections of $U_X$ through  $u\in U_X$ and their projections on $X$ (see the proof of Lemma 2 above), since  the statements are trivial in dimension one.
The required Hermitian quasi-metric on $\Cal L^q_X$ is also defined by 
induction on dimension with the help of the Poincar\'e residue map \cite{3,
pp.\;147-148}.

The condition $\ell^2$ on orbits of $\Gamma$ is a local property on $X$. We
get only a Hermitian quasi-metric on $\Cal L^q_X$ (instead of a Hermitian metric). Precisely, 
we get Hermitian metrics over small neighborhoods of  points of $X$, and on the intersections of neighborhoods, they will differ by constant multipliers (see \cite{5, Chap.\;5.13}).

For $\forall k>N\gg 0$, the Poincar\'e series are continuous sections
$$
P(f^k)(u):= \sum_{\gamma\in \Gamma} \gamma^*f^k(\gamma u),
$$
and they do not vanish for infinitely many $k$ (see \cite{5, Sect.\;13.1, Theorem 13.2}).

Finally, we can apply Gromov's theorem, precisely, its generalization by Koll\'ar 
(see \cite{4, Corollary 3.2.B, Remark 3.2.B$'$} and 
\cite{5, Theorem 13.8, Corollary 13.8.2, Theorem 13.9, Theorem 13.10}).
So, $X$ is a Moishezon manifold.

The Lemma 3 and Theorem are established.

\enddemo

\remark{Remarks} i) The theorem of the present note provides an alternative proof of a conjecture of H. Wu  provided $\pi_1(X)$ is residually finite (see \cite{10}).

ii) A generalization of the theorem to singular spaces will appear elsewhere.


\endremark

\enddocument